\documentclass[11pt]{amsart}

\newcommand{\BC}{{\mathbf C}}
\newcommand{\BF}{{\mathbf F}}
\newcommand{\BP}{{\mathbf P}}
\newcommand{\BQ}{{\mathbf Q}}
\newcommand{\BR}{{\mathbf R}}
\newcommand{\BZ}{{\mathbf Z}}

\newcommand{\Aid}{{\mathfrak A}}
\newcommand{\Oid}{{\mathfrak O}}

\newcommand{\zbar}{\overline{z}}

\newcommand{\Aut}{\operatorname{Aut}}
\newcommand{\Cl}{\operatorname{Cl}}

\newcommand{\GL}{\operatorname{GL}}

\newcommand{\PGL}{\operatorname{PGL}}

\newcommand\ra{\rightarrow}

\newcommand\buff{\rule[-1ex]{0pt}{3.3ex}}

\newtheorem{theorem}{Theorem}
\newtheorem{lemma}[theorem]{Lemma}

\newtheorem{proposition}[theorem]{Proposition}

\theoremstyle{definition}

\theoremstyle{remark}

\newtheorem{example}{Example}		
\newtheorem{acks}{Acknowledgments}	

\begin{document}

\title[Jacobians of plane quartics]
{Plane quartics with Jacobians isomorphic to a hyperelliptic Jacobian}
\subjclass{Primary 14H40; Secondary 14H45}
\keywords{Curve, Jacobian, polarization, Torelli, quartic}
\author{Everett W.\ Howe}
\address{Center for Communications Research, 
         4320 Westerra Court, 
         San Diego, CA 92121-1967, USA.}
\email{however@alumni.caltech.edu}


\date{27 May 1998}

\begin{abstract}
We show how for every integer $n$ one can explicitly construct $n$ distinct
plane quartics and one hyperelliptic curve over $\BC$ all of whose Jacobians
are isomorphic to one another as abelian varieties without polarization.  
When we say that the curves can be constructed ``explicitly'',
we mean that the coefficients of the defining equations of the curves
are simple rational expressions in algebraic numbers in $\BR$ whose minimal polynomials
over $\BQ$ can be given exactly
and whose decimal approximations can be given to as many places
as is necessary to distinguish them from their conjugates.
We also prove a simply-stated theorem that allows one to decide whether or not
two plane quartics over $\BC$, each with a pair of commuting involutions,
are isomorphic to one another.
\end{abstract}

\maketitle

\section{Introduction}

Torelli's theorem states that a curve is determined by its polarized Jacobian variety,
but a century ago Humbert showed that distinct curves can have isomorphic
unpolarized Jacobian varieties; 
thus it is natural to wonder exactly how much information about a curve is contained
in its unpolarized Jacobian.
In this paper we will prove that in general one cannot determine whether or not a curve
over the complex numbers is hyperelliptic simply by looking at its unpolarized Jacobian.
To prove this, we will show how for every positive integer $n$ one can
explicitly construct $n$ distinct plane quartics and one hyperelliptic curve
of genus~$3$ such that all $n+1$ of these curves share the same unpolarized Jacobian. 
Our construction is apparently the first method of producing explicit exact
equations for curves of genus~$3$ over~$\BC$ with isomorphic Jacobians.
In order to state our theorem precisely, we must set some notation.

For every nonzero complex number $\alpha$, we let $H(\alpha)$ denote the 
normalization of the possibly-reducible curve
\begin{align*}
W^2 &= X^4 + Y^4 + \left(-1 + \frac{1}{2\alpha^2}\right)\\
1   &= X^2 + Y^2.
\end{align*}
It is not hard to check that if $\alpha^2\neq 1$ then $H(\alpha)$ 
is irreducible and is in fact 
a hyperelliptic curve of genus~$3$.
For every pair of complex numbers $(\alpha,\beta)$ with
$\alpha\neq 0$ and $\alpha^2\neq 1$, we let $C(\alpha,\beta)$
denote the possibly-singular plane quartic curve defined by the homogeneous equation
$$X^4 + Y^4 + Z^4 + \left(-2 + 4\frac{1-\beta^2}{1-\alpha^2}\right) X^2Y^2 
   + \left(\frac{2\beta}{\alpha}\right) X^2Z^2 + \left(\frac{2\beta}{\alpha}\right) Y^2Z^2 = 0.$$
Finally, for every nonzero real $x$ we let $\mu(x)$ denote the unique $\mu\in\BR$ such that
$0\le \mu<1$ and such that the elliptic curve
$Y^2 = (X-1)(X-\mu)(X+1)$
is complex-analytically isomorphic to the torus $\BC/(\BZ + ix\BZ)$.

\begin{theorem}
\label{main}
Let $m$ be a positive even squarefree integer, and for every odd positive divisor $d$ of $m$
let $\alpha_d = \mu(\sqrt{m}/d)$.  For every odd divisor $d>1$ of $m$ 
the quartic $C(\alpha_1,\alpha_d)$ is nonsingular, and 
its Jacobian is isomorphic to the Jacobian of $H(\alpha_1)$
as an abelian variety without polarization.
Furthermore, if $d$ and $d'$ are distinct odd divisors of $m$ that
are greater than~$1$, then $C(\alpha_1,\alpha_d)$
and $C(\alpha_1,\alpha_{d'})$ are not isomorphic to one another.
\end{theorem}

When $x^2$ is a rational number the number $\mu(x)$ is algebraic and its minimal polynomial
over $\BQ$ can be calculated.
For instance, the values of $\mu(x)$ that we must calculate in order to
apply the theorem with $m=30$ are
\begin{align*}
\scriptscriptstyle  \mu(\sqrt{30})     &\scriptscriptstyle \ =\ - 464325           - 328320 \sqrt{2}  + 268072 \sqrt{3}
                                                 - 207648 \sqrt{5}  + 189560 \sqrt{6}  - 146832 \sqrt{10}
                                                 + 119888 \sqrt{15} +  84772 \sqrt{30}\\
\scriptscriptstyle  \mu(\sqrt{30}/3)   &\scriptscriptstyle \ =\  - 464325           + 328320 \sqrt{2}  - 268072 \sqrt{3}
                                                 + 207648 \sqrt{5}  + 189560 \sqrt{6}  - 146832 \sqrt{10}
                                                 + 119888 \sqrt{15} -  84772 \sqrt{30}\\
\scriptscriptstyle  \mu(\sqrt{30}/5)   &\scriptscriptstyle \ =\  - 464325           + 328320 \sqrt{2}  - 268072 \sqrt{3}
                                                 - 207648 \sqrt{5}  + 189560 \sqrt{6}  + 146832 \sqrt{10}
                                                 - 119888 \sqrt{15} +  84772 \sqrt{30}\\
\scriptscriptstyle -\mu(\sqrt{30}/15)  &\scriptscriptstyle \ =\  - 464325           - 328320 \sqrt{2}  + 268072 \sqrt{3}
                                                 + 207648 \sqrt{5}  + 189560 \sqrt{6}  + 146832 \sqrt{10}
                                                - 119888 \sqrt{15} -  84772 \sqrt{30}.
\end{align*}
Using these values, we find three distinct explicitly-given plane quartics whose
Jacobians are all isomorphic to that of an explicitly-given hyperelliptic curve.

For our proof of Theorem~\ref{main} we will require a simple method of determining
whether two plane quartics, each with a pair of commuting involutions, are
isomorphic to one another.  
We will provide such a method in Section~\ref{section-detecting}.
We will then prove Theorem~\ref{main} in Sections~\ref{section-isomorphic}
and~\ref{section-distinct}.
Finally, in Section~\ref{section-calculating} we will show how to
compute the minimal polynomials of the algebraic numbers $\mu(\sqrt{m}/d)$
that appear in the statement of Theorem~\ref{main}.
We will work through the calculations for the cases $m=6$ and $m=30$.

There are a number of other papers that discuss the relationship between a curve
and its unpolarized Jacobian.
The fact that a curve is not determined by its unpolarized Jacobian was first observed 
by Humbert~\cite{humbert}, who exhibited the period matrices for pairs of genus-$2$ curves
over~$\BC$ with isomorphic Jacobians.  
The Jacobians of Humbert's curves are reducible, and
Hayashida and Nishi \cite{hayashida-nishi, hayashida} 
showed that in fact there exist arbitrarily large sets of genus-$2$ curves over~$\BC$
all sharing the same reducible Jacobian.
A method for producing explicit equations for the curves in such sets was given in~\cite{howe:complex}.
There also exist {\it simple\/} abelian varieties over $\BC$ of dimension $2$, $3$, and~$4$
that can be obtained in more than one way as the
Jacobian of a curve~---~see Lange \cite{lange} for dimensions $2$ and $3$, 
and Ciliberto and van der Geer \cite{ciliberto-vandergeer} for dimension $4$.
No explicit examples of the equations for such curves over $\BC$ are known.
Over finite fields, explicit examples of distinct curves of genus~$2$ and~$3$ sharing the
same reducible Jacobian can be obtained by 
using results of Ibukiyama, Katsura, and Oort~\cite{ibukiyama-katsura-oort}
and Brock~\cite{brock}
or by reducing the examples in~\cite{howe:complex},
and a method for producing explicit examples of distinct curves
of genus~$2$ and~$3$ sharing the same irreducible Jacobian is given in~\cite{howe:finite}.
The genus-$3$ example worked out in~\cite{howe:finite} consists of a hyperelliptic
curve and a plane quartic over $\BF_3$ with isomorphic Jacobians.

\begin{acks}
The author thanks I.~M.~Isaacs for providing succinct descriptions for some of the
groups appearing in the proof of Proposition~\ref{curves}.
\end{acks}

\section{Detecting isomorphisms between plane quartics with commuting involutions}
\label{section-detecting}

Fix a set of homogeneous coordinates $X$, $Y$, $Z$ for $\BP^2$, and for every triple
$(a,b,c)$ of complex numbers let $Q(a,b,c)$ denote plane quartic defined by
$$X^4 + Y^4 + Z^4 + a X^2 Y^2 + b X^2 Z^2 + c Y^2 Z^2 = 0.$$
The curve $Q(a,b,c)$ is nonsingular if and only if 
$a^2 + b^2 + c^2 - abc - 4$ is nonzero and none of $a^2$, $b^2$, and $c^2$ is 
equal to $4$.  Our goal in this section is to give a simple criterion for deciding
whether two nonsingular curves $Q(a,b,c)$ and $Q(a',b',c')$ are isomorphic to one another.

The embedding of $Q(a,b,c)$ into $\BP^2$ given by its defining equation is a canonical embedding,
so every isomorphism $\varphi$ from $Q(a,b,c)$ to $Q(a',b',c')$ can be extended to give an
automorphism $\varphi_{\BP^2}$ of the ambient $\BP^2$ that takes $Q(a,b,c)$ to $Q(a',b',c')$.
Using our fixed set of homogeneous coordinates, we can identify $\Aut \BP^2$ with $\PGL(3,\BC)$,
so $\varphi_{\BP^2}$ can be represented by a $3\times 3$ matrix, unique up to scalar multiples.
We say that an isomorphism $\varphi\colon Q(a,b,c)\ra Q(a',b',c')$ is {\it strict\/}
if $\varphi_{\BP^2}$ has a representative that is the product of a permutation matrix and a diagonal matrix.
It is easy to see that $Q(a',b',c')$ is strictly isomorphic to $Q(a,b,c)$ if and only if
the triple $(a',b',c')$ can be obtained from $(a,b,c)$ by permuting the order of
the elements and changing the signs of an even number of elements.

\begin{proposition}
\label{curves}
An isomorphism class of nonsingular quartics of the form $Q(a',b',c')$ is equal to one of the
following{\rm:}
\begin{enumerate}
\item[1.] the strict isomorphism class of the curve $Q(a,b,c)$ for some $a$, $b$, and~$c$
such that $a^2$, $b^2$, and~$c^2$ are pairwise unequal{\rm;}
\item[2.] the union of the strict isomorphism classes of the curves
$Q(a,b,b)$ and $Q(-2 + 16/(a+2), 2b/d,2b/d)$
for some $a$ and $b$ with $b\neq 0$, where $d^2 = a+2${\rm;} or
\item[3.] the union of the strict isomorphism classes of the curves
$Q(a,0,0)$ and $Q( -2 + 16/(a+2), 0,0,)$ and
$Q( -2 + 16/(-a+2), 0,0,)$, for some $a$.
\end{enumerate}
\end{proposition}

A special case of this proposition may be found in Kuribayashi and Sekita~\cite{kuribayashi-sekita},
but we have been unable to find a proof of the general case in the literature.
Brock uses the result of the proposition, without proof, in Chapter~3 of~\cite{brock}.

Our proof of the proposition depends on the following lemma.
By a {\it $V_4$-subgroup\/} of a group $G$, we mean a subgroup of $G$
isomorphic to the Klein $4$-group $V_4$. 

\begin{lemma}
\label{V4-subgroups}
Let $C$ be a nonhyperelliptic curve of genus~$3$.  The number of strict isomorphism
classes of curves $Q(a,b,c)$ that are isomorphic to $C$ is equal to the number
of conjugacy classes of $V_4$-subgroups of $\Aut C$.
\end{lemma}

\begin{proof}
For every nonsingular curve $Q(a,b,c)$
let $A(a,b,c)$ be the subgroup of $\Aut Q(a,b,c)$ consisting of
the automorphisms $[X\!:\! Y\!:\! Z]\mapsto[\pm X\!:\! \pm Y\!:\! \pm Z]$, so that $A(a,b,c)$ is
isomorphic to $V_4$.
If $\varphi\colon C\ra Q(a,b,c)$ is an isomorphism, then $\varphi^*A(a,b,c)$
is a $V_4$-subgroup of $\Aut C$.
If $\psi$ is another isomorphism from $C$ to the same curve $Q(a,b,c)$, then $\psi^*A(a,b,c)$
and $\varphi^*A(a,b,c)$ are conjugate subgroups of $\Aut C$.
Furthermore, if $\psi\colon C\ra Q(a',b',c')$ is the composition of $\varphi$ with
a strict isomorphism from $Q(a,b,c)$ to $Q(a',b',c')$, then $\varphi^*A(a,b,c) = \psi^*A(a',b',c')$.  
Thus we get a map $\Theta$ from the set of strict isomorphism classes of curves $Q(a,b,c)$
such that $C\cong Q(a,b,c)$ to the set of conjugacy classes of $V_4$-subgroups of $\Aut C$.
We will prove that $\Theta$ is a bijection.

The canonical embedding $\Phi\colon C\ra \BP^2$ allows us to identify
$\Aut C$ with a subgroup of $\Aut \BP^2$.
Suppose $A$ is a $V_4$-subgroup of $\Aut C$.
Since there is a unique embedding of $V_4$ into $\Aut \BP^2$
up to conjugacy, we may choose coordinates $X$, $Y$, $Z$ for $\BP^2$
so that the image of $A$ in $\Aut \BP^2$ consists of the automorphisms
$[X\!:\! Y\!:\! Z]\mapsto[\pm X\!:\! \pm Y\!:\! \pm Z]$.
It is easy to see that with this choice of
coordinates for~$\BP^2$, the image of $C$ under the canonical embedding is a plane
quartic defined by a homogeneous quadratic polynomial in $X^2$, $Y^2$, and $Z^2$.
By rescaling coordinates, we find that $\Phi(C)$ is defined by a quartic of the
form $Q(a,b,c)$ for some $a$, $b$, and $c$. By construction, the
isomorphism $\varphi\colon C\ra Q(a,b,c)$ obtained in this way pulls back $A(a,b,c)$
to our original group $A$. Thus $\Theta$ is surjective.

On the other hand,
suppose $\varphi\colon C\ra Q(a,b,c)$ and $\psi\colon C\ra Q(a',b',c')$ are
isomorphisms such that $\varphi^* A(a,b,c)$ is conjugate to $\psi^* A(a',b',c')$,
say by an element $\alpha\in\Aut C$.
By replacing $\psi$ with the composition $\psi\alpha$, we may assume
that $\varphi^* A(a,b,c)$ is equal to $\psi^* A(a',b',c')$.
Let $\chi\colon Q(a,b,c)\ra Q(a',b',c')$ be the isomorphism $\psi\varphi^{-1}$.
Then the fact that $\chi^*A(a',b',c') = A(a,b,c)$ shows that the
element $\chi_{\BP^2}$ of $\PGL(3,\BC)$ can be represented by
the product of a permutation matrix and a diagonal matrix.
Thus $Q(a',b',c')$ is strictly isomorphic to $Q(a,b,c)$,
so $\Theta$ is injective.
\end{proof}

\begin{proof}[Proof of Proposition~{\rm{\ref{curves}}}]
There are exactly seven groups containing $V_4$ that occur as automorphism groups of 
nonhyperelliptic curves of genus~$3$ over $\BC$ ---
see~\cite{kuribayashi-komiya-2}, or Theorem~5.5 of~\cite{vermeulen}, or Theorem~3.5 of~\cite{brock}.
These groups are $V_4$ itself, the dihedral group
$D_8$ of order $8$, the symmetric group $S_4$, 
a group $G_{16}$ that is isomorphic to the 
central product of the quaternion group $Q_8$ with the cyclic group $C_4$
(that is, the quotient of the product $Q_8\times C_4$ 
by the image of a central diagonal embedding of $C_2$), 
a group $G_{48}$ that is isomorphic to the semidirect product of $G_{16}$ with $C_3$
(where the $C_3$ acts in the obvious way on the quaternion group),
a group $G_{96}$ that is isomorphic to the semidirect product of 
the trace-$0$ part of $(\BZ/4\BZ)\times(\BZ/4\BZ)\times(\BZ/4\BZ)$
with $S_3$ (where $S_3$ acts by permuting the factors),
and the simple group $\GL(3,\BF_2)$ of order $168$.

For each of these groups $G$ we can calculate the number $N$
of conjugacy classes of $V_4$-subgroups of $G$.
We leave these straightforward calculations to the reader; the results are presented in the first
two columns of Table~\ref{table}.
For every automorphism group $G$, Vermeulen (\cite{vermeulen}, Table~5.6, pp.~63--64)
lists a standard way of writing the curves with that automorphism group.  
For all groups in Table~\ref{table} except for $D_8$ and $\GL(3,\BF_2)$
we list Vermeulen's standard form in column~$3$.  For $D_8$ we list a standard form
easily obtained from Vermeulen's, and for $\GL(3,\BF_2)$ we list the form 
that was apparently first obtained by Ciani \cite{ciani}.
Now, if a curve $Q(a,b,c)$ has automorphism group equal to $G$, then 
Lemma~\ref{V4-subgroups} says there there will be exactly $N$
strict isomorphism classes of curves $Q(a',b',c')$ isomorphic to $Q(a,b,c)$.
One of these classes will be represented by $Q(a,b,c)$ itself.  In
column~$4$ we list representatives of the other $N-1$ strict isomorphism classes
of curves $Q(a',b',c')$ isomorphic to $Q(a,b,c)$,
all of which are obtained by applying Lemma~\ref{vary-the-curve} (below) to elements 
of the strict isomorphism class of $Q(a,b,c)$.
The reader may verify that the strict isomorphism classes of these $N$ curves are distinct
from one another whenever $Q(a,b,c)$ has automorphism group exactly $G$.

Proposition~\ref{curves} follows immediately upon inspection of Table~\ref{table}.
\end{proof}

\begin{table}
\begin{center}
\begin{tabular}{|c|c|c|c|}
\hline
$G$            & $N$ & Standard form        & Associated form(s)           \buff \\ \hline\hline
$V_4$          & $1$ & $Q(a,b,c)$           & none                         \buff \\ \hline
$D_8$          & $2$ & $Q(a,b,b)$           & $Q(-2+16/(a+2), 2b/d, 2b/d)$ \buff \\
               &     &                      & (where $d^2 = a+2$)          \buff \\ \hline
$S_4$          & $2$ & $Q(a,a,a)$           & $Q(-2+16/(a+2), 2a/d, 2a/d)$ \buff \\ 
               &     &                      & (where $d^2 = a+2$)          \buff \\ \hline
$G_{16}$       & $3$ & $Q(a,0,0)$           & $Q(-2+16/(a+2),0,0)$         \buff \\
               &     &                      & $Q(-2+16/(-a+2),0,0)$        \buff \\ \hline
$G_{48}$       & $1$ & $Q(2\sqrt{-3},0,0)$  & none                         \buff \\ \hline
$G_{96}$       & $2$ & $Q(0,0,0)$           & $Q(6,0,0)$                   \buff \\ \hline
$\GL(3,\BF_2)$ & $2$ & $Q(z,z,z)$           & $Q(\zbar,\zbar,\zbar)$       \buff \\ \hline
\end{tabular}
\end{center}
\vspace{1ex}
\caption{Forms of curves with automorphism groups containing commuting involutions.
For every possible automorphism group $G$ that contains a $V_4$-subgroup
we list the number $N$ of conjugacy classes of such subgroups. 
Every curve $C$ with $G\subseteq \Aut C$ can be put in the standard form
$Q(\cdot,\cdot,\cdot)$ listed in the third column.  
If $Q(a,b,c)$ has automorphism group $G$, then the $N-1$ other strict isomorphism
classes of curves $Q(a',b',c')$ isomorphic to $Q(a,b,c)$ are listed in column~$4$.
The number $z$ in the last row is $z = 3(-1 + \protect\sqrt{-7})/2$, and $\zbar$ 
is its complex conjugate.}
\label{table}
\end{table}


\begin{lemma}
\label{vary-the-curve}
Suppose $Q(a,b,b)$ is a nonsingular curve. Let $d\in \BC$ satisfy $d^2 = a+2$.
Then $Q(-2+16/(a+2),2b/d,2b/d)$ is isomorphic to $Q(a,b,b)$.
\end{lemma}

\begin{proof}
Let $e\in\BC$ satisfy $e^2 = d$.
Then the element of $\Aut \BP^2 \cong\PGL(3,\BC)$ represented by the matrix
$$\left(
\begin{matrix}
e/2 & e/2  & 0 \\
e/2 & -e/2 & 0 \\
0   & 0    & 1
\end{matrix}
\right)$$
gives an isomorphism from $Q(a,b,b)$ to $Q(-2+16/(a+2),2b/d,2b/d)$.
\end{proof}

\section{Proof that the curves in Theorem~\ref{main} have isomorphic Jacobians}
\label{section-isomorphic}

We will present three lemmas that together will prove that 
$C(\alpha_1,\alpha_d)$ is nonsingular
and that $C(\alpha_1,\alpha_d)$ and $H(\alpha_1)$
have isomorphic Jacobians, for all odd divisors $d>1$ of $m$.
The proofs of the first two lemmas rely on the results of Part~4 of~\cite{howe-leprevost-poonen}, which tell us
how to find genus-$3$ curves whose Jacobians are isomorphic to abelian varieties
of the form $(E_1\times E_2\times E_3)/G$ for elliptic curves $E_i$ and
for subgroups $G$ of a certain type.  

For every odd positive divisor $d$ of $m$ let $F_d$ be the 
elliptic curve given by $Y^2 = X(X^2 + 2\alpha_d X + \alpha_d^2 -1)$.
Also, let $\beta = 2\alpha_1^2 - 1$ and let $F'$ be the elliptic curve 
$Y^2 = X(X^2 + 2\beta X + \beta^2 - 1)$.
Let $S_d$, $T_d$, and $U_d$ be the $2$-torsion points on $F_d(\BC)$
with $x$-coordinates $-1-\alpha_d$,~$0$, and~$1-\alpha_d$, respectively, and let
$S'$, $T'$, and $U'$ be the $2$-torsion points on $F'(\BC)$ with $x$-coordinates
$-1-\beta$,~$0$, and~$1-\beta$, respectively.
For every $d$, let $G_d$ be the subgroup of $(F_d\times F_d\times F')(\BC)$
generated by $(T_d,0,T')$, $(0,T_d,T')$, and~$(S_d, S_d, S')$.

\begin{lemma} 
\label{hyperelliptic}
The Jacobian of $H(\alpha_1)$ is isomorphic to $(F_1\times F_1\times F')/G_1$.
\end{lemma}

\begin{proof}
We will find a curve whose Jacobian is isomorphic to $(F_1\times F_1\times F')/G_1$
by applying the construction of Section~4.1 of~\cite{howe-leprevost-poonen},
which requires that we specify three elliptic curves $E_1$, $E_2$, $E_3$ and
$2$-torsion points $P_i$ and $Q_i$ on each $E_i$.
We take $(E_1, P_1, Q_1)$ to be $(F_1, S_1, T_1)$, we take
$(E_2, P_2, Q_2)$ to be $(F_1, S_1, T_1)$, and we take
$(E_3, P_3, Q_3)$ to be $(F', S', T')$.
Note that then the group $G$ of~\cite{howe-leprevost-poonen}
is equal to our group $G_1$.

In the notation of Part~4 of~\cite{howe-leprevost-poonen} we have 
$$\begin{matrix}
A_1 = 2\alpha_1,\ & B_1 = \alpha_1^2 - 1,\ & 
   \Delta_1 = 4,\ & d_1 = 2,\ & \text{and\ } \lambda_1 = \alpha_1,\buff\\
A_2 = 2\alpha_1,\ & B_2 = \alpha_1^2 - 1,\ & 
   \Delta_2 = 4,\ & d_2 = 2,\ & \text{and\ } \lambda_2 = \alpha_1,\buff\\
A_3 = 2\beta,\    & B_3 = \beta^2 - 1,\    & 
   \Delta_3 = 4,\ & d_3 = 2,\ & \text{and\ } \lambda_3 = \beta,\buff
\end{matrix}$$
we have $R=8$, and the ``twisting factor'' $T$ is $0$.  We can therefore apply Proposition~14 
of~\cite{howe-leprevost-poonen}, and we find that 
$(F_1\times F_1\times F')/G_1$ is isomorphic to the Jacobian of the curve
defined (in homogeneous coordinates) by 
\begin{align*}
W^2 Z^2 & = a X^4 + b Y^4 + c Z^4\\
      0 & = d X^2 + e Y^2 + f Z^2
\end{align*}
where 
$$\begin{matrix}
\buff a = 4\alpha_1^2(\alpha_1^2-1)^2\hfill                     & d = 1/(2\alpha_1(\alpha_1^2-1))\hfill \\
\buff b = 4\alpha_1^2(\alpha_1^2-1)^2\hfill                     & e = 1/(2\alpha_1(\alpha_1^2-1))\hfill \\
\buff c = -8\alpha_1^2(\alpha_1^2-1)^2(2\alpha_1^2-1)\ \ \hfill & f = -1/(\alpha_1^2-1).\hfill
\end{matrix}$$
If we divide the quartic equation by $a$ and the quadratic by $d$,
replace $Z$ with $Z/\sqrt{2\alpha_1}$, 
replace $W$ with $W\cdot 2\alpha_1(1-\alpha_1^2)\sqrt{2\alpha_1}$,
and dehomogenize with respect to $Z$, 
we find that the curve given to us by  Proposition~14 
of~\cite{howe-leprevost-poonen} is none other than $H(\alpha_1)$.
\end{proof}

\begin{lemma} 
\label{quartic}
If $d$ is a divisor of $m$ with $d>1$, then the quartic $C(\alpha_1,\alpha_d)$ is nonsingular
and its Jacobian is isomorphic to $(F_d\times F_d\times F')/G_d$.
\end{lemma}

\begin{proof}
We will again 
use the construction of Section~4.1 of~\cite{howe-leprevost-poonen}.
This time we take $(E_1, P_1, Q_1)$ to be $(F_d, S_d, T_d)$, we take
$(E_2, P_2, Q_2)$ to be $(F_d, S_d, T_d)$, and we take
$(E_3, P_3, Q_3)$ to be $(F', S', T')$.
Now the group $G$ of~\cite{howe-leprevost-poonen}
is equal to our group $G_d$.
In the notation of Part~4 of~\cite{howe-leprevost-poonen} we have
$$\begin{matrix}
\buff A_1 = 2\alpha_d,\ & B_1 = \alpha_d^2 - 1,\ & 
   \Delta_1 = 4,\ & d_1 = 2,\ & \text{and\ } \lambda_1 = \alpha_d,\\
\buff A_2 = 2\alpha_d,\ & B_2 = \alpha_d^2 - 1,\ & 
   \Delta_2 = 4,\ & d_2 = 2,\ & \text{and\ } \lambda_2 = \alpha_d,\\
\buff A_3 = 2\beta,\    & B_3 = \beta^2 - 1,\    & 
   \Delta_3 = 4,\ & d_3 = 2,\ & \text{and\ } \lambda_3 = \beta,
\end{matrix}$$
we have $R=8$, and the twisting factor $T$ is $8(\beta-1)(\beta-2\alpha_d^2+1)$.
We see that $T$ is nonzero, because $\beta = 2\alpha_1^2 - 1 <1$
and because $\alpha_1\neq\pm\alpha_d$.
We can therefore apply Proposition~15
of~\cite{howe-leprevost-poonen} to find that 
$(F_d\times F_d\times F')/G_d$ is isomorphic to the Jacobian of the
nonsingular plane quartic
$$B_1 X^4 + B_2 Y^4 + B_3 Z^4 + dX^2Y^2 + eX^2Z^2 + fY^2Z^2 = 0$$
where
$$ d = 2 (1-\alpha_d^2) - 4(1-\alpha_1^2)
\text{\qquad and\qquad} e  =  f  = 4\alpha_d(1-\alpha_1^2).$$
If we multiply the equation for the quartic by $-1$, 
replace $X$ with $X/(1-\alpha_d^2)^{1/4}$, 
replace $Y$ with $Y/(1-\alpha_d^2)^{1/4}$, 
and replace $Z$ with $Z/(1-\beta^2)^{1/4}$, 
we find that this nonsingular plane quartic is isomorphic to the curve given by 
$$X^4 + Y^4 + Z^4 + d' X^2 Y^2 + e' X^2 Z^2 + e' Y^2 Z^2 = 0,$$
where
$$d' = -2 + 4\frac{1-\alpha_1^2}{1-\alpha_d^2}
\text{\qquad and \qquad}
e' = -2\frac{\alpha_d}{\alpha_1} \sqrt{\frac{1-\alpha_1^2}{1-\alpha_d^2}}.$$
Finally, by applying Lemma~\ref{vary-the-curve} 
we find that this last curve is isomorphic to $C(\alpha_1,\alpha_d)$.
\end{proof}

\begin{lemma}
\label{isomorphic}
For every odd positive divisor $d$ of $m$ we have
$$(F_1\times F_1\times F')/G_1 \cong (F_d\times F_d\times F')/G_d$$
as abelian varieties without polarization.
\end{lemma}

\begin{proof}
By shifting $X$-coordinates by $\alpha_1$, we see that $F_1$ is
isomorphic to the curve $Y^2 = (X-1)(X-\alpha_1)(X+1)$.
Since $\alpha_1 = \mu(\sqrt{m})$, the
definition of $\mu$ shows that $F_1$ is complex-analytically
isomorphic to the torus $\BC/\Lambda_1$, where
$\Lambda_1 = \BZ + i\sqrt{m}\BZ$.
This isomorphism is given by sending a point $P$ on $F_1$ to the image in $\BZ/\Lambda_1$
of the integral from $\infty$ to $P$ of $k\ dX/Y$ for some constant $k$ that is either 
real or pure imaginary.
Since $1/Y$ is pure imaginary for real values of $X$ less than $-1-\alpha_1$
and real for real values of $X$ between $-1-\alpha_1$ and~$0$,
we see that $S_1$ corresponds to either $1/2+\Lambda_1$ or $i\sqrt{m}/2 + \Lambda_1$
and that $T_1$ corresponds to $(1 + i\sqrt{m})/2 + \Lambda_1$.
Similarly, if we let $\Lambda_d$ be the lattice $\BZ + i\sqrt{m}/d\BZ$, 
there is a complex-analytic isomorphism from $F_d$ to
$\BC/\Lambda_d$ that takes $S_d$ to either 
$1/2+\Lambda_d$ or $i\sqrt{m}/2d + \Lambda_d$
and that takes $T_d$ to $(1 + i\sqrt{m}/d)/2 + \Lambda_d$.

Let $u$ and $v$ be integers such that  $ud + v(m/d) = 1$  and such that 
$v$ is a multiple of $4$.  One can easily check that the matrix
$$\left(
\begin{matrix}
ud & i2\sqrt{m}\\
i(v/2)\sqrt{m} & d
\end{matrix}
\right)$$
gives an automorphism of $\BC\times\BC$ that takes $\Lambda_1\times\Lambda_1$
to $\Lambda_d\times\Lambda_d$.
Thus, this matrix gives an isomorphism from $(\BC/\Lambda_1)\times(\BC/\Lambda_1)$
to $(\BC/\Lambda_d)\times(\BC/\Lambda_d)$,
which we may interpret as an isomorphism $\varphi\colon F_1\times F_1\ra F_d\times F_d$.
Furthermore, using the matrix interpretation of $\varphi$, it is easy to check that
$\varphi((T_1,0)) = (T_d,0)$, that
$\varphi((0,T_1)) = (0,T_d)$, and that
$\varphi((S_1,S_1))$ is either $(S_d,S_d)$ or $(U_d,U_d)$.
But then $\varphi\times 1_{F'}$ is an isomorphism from $F_1\times F_1\times F'$
to $F_d\times F_d\times F'$ that takes $G_1$ to~$G_d$, and we are done.
\end{proof}

\section{Proof that the plane quartics in Theorem~\ref{main} are distinct}
\label{section-distinct}

To prove that distinct values of $d$ give us distinct curves $C(\alpha_1, \alpha_d)$, we will
need the following lemma.

\begin{lemma}
\label{one-largest}
For every odd divisor $d>1$ of $m$ we have $\alpha_1 > \alpha_d$.
\end{lemma}

\begin{proof}
The $j$-invariant of the elliptic curve $Y^2 = (X-1)(X-\mu)(X+1)$
is $2^6(\mu^2+3)^3/(\mu^2-1)^2$,
and this is an increasing function of $\mu$ for $\mu\in[0,1)$.
Likewise, the function from $\BR$ to $\BR$ that takes
$x$ to the $j$-invariant of the lattice $\BZ + ix\BZ$ is increasing
for $x\in[1,\infty)$.
Thus, the function $\mu\colon\BR\ra\BR$ defined in the introduction
is increasing for $x\in[1,\infty)$.
Suppose $d>1$ is an odd divisor of $m$, so that $1 < d < m$.
If $d\le\sqrt{m}$ then $1\le \sqrt{m}/d < \sqrt{m}$
and we have
$\mu(\sqrt{m}/d) < \mu(\sqrt{m})$.
If $d > \sqrt{m}$ then $1 < d/\sqrt{m} < \sqrt{m}$,
and by using the fact that $\mu(x) = \mu(1/x)$ we find that 
$\mu(\sqrt{m}/d) = \mu(d/\sqrt{m}) < \mu(\sqrt{m})$.
\end{proof}

Note that $\alpha_d\neq 0$ for every $d$, because $\mu(x)=0$ only for $x=1$. From
Proposition~\ref{curves} we see
that if two curves $Q(a,b,b)$ and $Q(a',b',b')$ 
are isomorphic to one another and if $bb'\neq0$,
then either $a=a'$ or $(a+2)(a'+2) = 16$.
So if our curves $C(\alpha_1,\alpha_d)$ and $C(\alpha_1,\alpha_{d'})$
were isomorphic to one another, we would find that either 
$$\frac{1-\alpha_1^2}{1-\alpha_d^2}=\frac{1-\alpha_1^2}{1-\alpha_{d'}^2}
\text{\qquad or \qquad}
\left(\frac{1-\alpha_1^2}{1-\alpha_d^2}\right)\left(\frac{1-\alpha_1^2}{1-\alpha_{d'}^2}\right) = 1.$$
The first equality is impossible because $\alpha_d$ and $\alpha_{d'}$
are distinct positive numbers and $\alpha_1^2\neq 1$, and the second equality is
impossible because Lemma~\ref{one-largest} shows that 
the two factors on the left-hand side
are both positive numbers less than $1$.
Thus the two curves are not isomorphic to one another, and
we have proven Theorem~\ref{main}.

\section{Explicitly calculating the equations for the curves}
\label{section-calculating}

Let $m$ be a positive even squarefree integer, let $K$ be the field $\BQ(\sqrt{-m})$,
and let $\Oid = \BZ[\sqrt{-m}]$ be the ring of integers of $K$. 
By choosing one of the two embeddings $K\hookrightarrow \BC$ we can think
of $K$ and $\Oid$ as subrings of $\BC$.  From the explicit class field theory
of imaginary quadratic extensions of $\BQ$, we know that there is a monic irreducible polynomial
$f\in\BZ[x]$ whose roots are the $j$-invariants of the elliptic curves over $\BC$
with complex multiplication by $\Oid$.  These elliptic curves correspond to the lattices in
$K\subset\BC$ given by fractional $\Oid$-ideals, with two lattices giving the same elliptic curve if
and only if the corresponding ideals represent the same element of the class group $\Cl K$
of~$K$.
Thus, by computing ideals that represent the various elements of the class group of $K$
and by computing to high-enough precision the $j$-invariants of the lattices given by these ideals,
we can find the polynomial $f$ {\it exactly\/} by first determining its coefficients
with an error of less than $1/2$ and then making use of the fact that these coefficients are integers.

For every odd positive divisor $d$ of $m$ let $\Aid_d$ be the ideal $(d,\sqrt{-m})$.
The ideals $\Aid_d$ represent distinct $2$-torsion elements of $\Cl K$, 
and every $2$-torsion element is represented by an $\Aid_d$.
The $j$-invariant of the lattice $\Aid_d$ is $j(\sqrt{-m}/d)$,
so the polynomial $f$ is the minimal polynomial of the conjugate algebraic
integers $j(\sqrt{-m}/d)\in\BR$.  In fact, the $j(\sqrt{-m}/d)$ are
precisely the real roots of $f$ --- see statement~5.4.3, p.~124, of \cite{shimura}.

Now, the $j$-invariant of the elliptic curve $Y^2 = (X-1)(X-\mu)(X+1)$ is
$2^6(\mu^2 + 3)^3/(\mu^2 - 1)^2,$
so the numbers $\alpha_d= \mu(\sqrt{m}/d)$ of Section~\ref{section-isomorphic}
are zeros of the rational function
$f\left(2^6(x^2+3)^3/(x^2-1)^2\right).$
Suppose we write this rational function as $g(x)/h(x)$, 
where $g$ and $h$ are coprime elements of $\BZ[x]$ and $h$ is monic.
Then $g(\alpha_d) = 0$ for each positive odd divisor $d$ of $m$,
and in fact the $\alpha_d$ are precisely the real roots of $g$ that lie
between $0$ and $1$.
Thus we can specify the $\alpha_d$ precisely, by specifying $g$ and by specifying
each $\alpha_d$ to enough decimal places to distinguish it from the other $\alpha_d$'s.

\begin{example}
Consider the case $m=6$ of the preceding discussion.
The class group of $K = \BQ(\sqrt{-6})$ has order $2$, and 
its elements are represented by $\Aid_1$ and $\Aid_3$.
One calculates (using PARI/GP, for example) that
$j(\Aid_1) \approx 4831907.903351340$ and that
$j(\Aid_3) \approx    3036.096648660$.
Thus $$f\approx x^2 - 4834944.0x + 14670139392.0$$ so that 
$f = x^2 - 4834944x + 14670139392$.  
We find that 
\begin{align*}
g & = 2^{12}(x^4 - 36x^2 + 36)(x^4 + 276x^3 + 342x^2 - 396x - 207)\\
  &\qquad\qquad\times(x^4 - 276x^3 + 342x^2 + 396x - 207).
\end{align*}
The roots of $g$ between $0$ and $1$ are roots of the last two factors, and we find that
\begin{align*}
\mu(\sqrt{6})    & = -69 - 48 \sqrt{2} + 40 \sqrt{3} + 28 \sqrt{6}  \approx \phantom{-}0.9854941\\
-\mu(\sqrt{6}/3) & = -69 + 48 \sqrt{2} + 40 \sqrt{3} - 28 \sqrt{6}  \approx           -0.4214295
\end{align*}
(where we list $-\mu(\sqrt{6}/3)$ instead of $\mu(\sqrt{6}/3)$ to emphasize the
fact that the two given values are conjugate algebraic integers).
\end{example}

\begin{example}
Now let us take $m$ to be $30$.
The class group of $K = \BQ(\sqrt{-30})$ has order $4$,
and its elements are represented
by the ideals $\Aid_1$, $\Aid_3$, $\Aid_5$, and $\Aid_{15}$.
One calculates that
\begin{align*}
j(\Aid_1)    & \approx            883067941288166.389365302091877037343676871900 \\
j(\Aid_3)    & \approx  \phantom{8830679412}96685.684561070142571546116321742465 \\
j(\Aid_5)    & \approx  \phantom{88306794128}1944.822048602018380008403472146653 \\
j(\Aid_{15}) & \approx  \phantom{8830679}29717203.104025025747171408136529238982.
\end{align*}
As before, we can find $f$ exactly and then compute $g$.  We find that
\begin{align*}
g &= 2^{24} (x^8 - 464328x^6 + 3576024x^4 - 6223392x^2 + 3111696)\\
  &\qquad  \times (x^8 - 3714600x^7 + 158287932x^6 - 550595160x^5 - 358871706x^4 \\
  &\qquad\qquad\qquad          + 1052916840x^3 + 302346108x^2 - 498574440x - 101729439)\\
  &\qquad \times   (x^8 + 3714600x^7 + 158287932x^6 + 550595160x^5 - 358871706x^4 \\
  &\qquad\qquad\qquad          - 1052916840x^3 + 302346108x^2 + 498574440x - 101729439)
\end{align*}
and the roots of $g$ between $0$ and $1$ are roots of the last two factors.
One can check that these roots are indeed the values of $\mu(\sqrt{30}/d)$ given in the introduction.
\end{example}

\end{document}